# A DESCENT THEOREM IN TOPOLOGICAL K-THEORY

## by Max Karoubi

**1. Introduction.** Let A be a real Banach algebra. We define its "K-theory space" $\mathcal{K}(A)$ as the direct limit

$$\operatorname*{colim}_{n} \operatorname{Proj}_n(A^{2n})$$

where $\operatorname{Proj}_n(A^{2n})$ denotes the space of projection operators in $A^{2n}$ with the 2n x 2n matrix

$$\begin{pmatrix} 1 & 0 \\ 0 & 0 \end{pmatrix}$$

as the base point. This space is homotopically equivalent (non canonically[1]) to the product

$$\mathcal{K}(A) = K_0(A) \times BGL(A)$$

where GL(A) and BGL(A) have the usual topology and $K_0(A)$ the discrete topology. Therefore, its homotopy groups are the topological K-groups of the Banach algebra A, which are periodic of period 8 (Bott). If A can be provided with a complex structure, these homotopy groups are in fact periodic of period 2.

**2.** If $A' = A \otimes_{\mathbf{R}} \mathbf{C}$ is the complexification of A, the group $G = \mathbf{Z}/2$ acts on $\mathcal{K}(A')$ by complex conjugation (we keep $G = \mathbf{Z}/2$ through all the paper) and we have a natural map

$$\sigma : \mathcal{K}(A) \longrightarrow \mathcal{K}(A')^{hG}$$

Here, $Y^{hG}$ denotes in general the homotopy fixed point set of the G-space Y. More precisely, $Y^{hG}$ is the space of (continuous) sections of the Borel fibration

$$\begin{array}{c} EG \times_G Y \\ \downarrow \\ BG \end{array}$$

It is easy to see that $Y^{hG}$ is also the space of equivariant maps $EG \longrightarrow Y$. The purpose of this paper is to prove the following "descent theorem"[2] :

---

[1] as it was pointed out to me by B. Kahn.

[2] This theorem is well known if A is the Banach algebra of real numbers : cf. [2], lemma 3.5 for instance. Other proofs have been given by J. Lannes (unpublished) and B. Kahn in a joint work with Hinda Hamraoui (in preparation). On the other hand, the statement seems new if A is the algebra of quaternions.

**3. THEOREM.** *The map* $\sigma : \mathcal{K}(A) \longrightarrow \mathcal{K}(A')^{hG}$ *defined above is a homotopy equivalence.*

**APPLICATIONS.**
One application of this theorem we have in mind is the Baum-Connes conjecture in the real case (this will be including in a forthcoming paper "BKR" by Paul Baum, John Roe and the author). The classical Baum-Connes conjecture (in the complex case) states that the index map

$$\mu(\Gamma) : K_j^{\Gamma}(\underline{E}\Gamma) \longrightarrow K_j(C_r^*(\Gamma))$$

is an isomorphism for a discrete (countable) group $\Gamma$ (where $j = 0, 1$ mod. 2). In the BKR paper, we show, using Theorem 3 among one of the essential ingredients, that if $\mu(\Gamma)$ is an isomorphism, the real analog

$$\mu_{\mathbf{R}}(\Gamma) : KO_j^{\Gamma}(\underline{E}\Gamma) \longrightarrow K_j(C_r^*(\Gamma, \mathbf{R}))$$

is also an isomorphism for $j = 0, 1..., 7$ mod. 8.

Another application is a comparison theorem between Algebraic K-theory and KR-theory of Real varieties (related to the Quillen-Lichtenbaum conjecture for such varieties). This will appear in a forthcoming paper by Charles Weibel and the author.

**4.** The first step of the proof of Theorem 3 is to remark the following basic fact : if Y is the product X x X with the action of G switching the factors, $Y^{hG}$ may be identified with the space of maps from EG to X, a space which is homotopically equivalent to X, since EG is contractible. Therefore, the map

$$\sigma : X \longrightarrow (X \times X)^{hG}$$

is a homotopy equivalence. From this remark, we deduce immediately the following proposition :

**5. PROPOSITION.** *The theorem is true if* A *is the underlying real algebra of a complex Banach algebra.*

*Proof.* It is well known and easy to see that the homomorphism $A' = A \otimes_{\mathbf{R}} \mathbf{C} \longrightarrow A \times A$ defined by $a \otimes z \mapsto (az, a\bar{z})$ is an isomorphism, the complex conjugation switching the factors. Therefore the space $\mathcal{K}(A')$ is just the product $\mathcal{K}(A) \times \mathcal{K}(A)$ and the proposition follows from the considerations in § 4.

The paper is now devoted to reduce the theorem to this easy case, using in an essential way the KR-theory of Atiyah [1]. The main ingredient in the proof is mentionned at the end of § 7.

**6. KR-theory of Banach algebras.** Let us consider a compact space X provided with an involution $x \mapsto \bar{x}$, following the notation used by Atiyah. We write A(X) for the Banach algebra of continuous functions $f : X \longrightarrow A'$ such that $f(\bar{x}) = \overline{f(x)}$, the complex conjugate of f(x). If A is the field of real numbers and X an arbitrary G-space, it is easy to see that the K-theory of A(X) is isomorphic to Atiyah's KR(X). If X is a locally compact space, we may extend this definition by choosing A(X) to be the space of continuous functions (with values in A') which go to 0 when x goes to ∞ (with the same conjugation condition). Note that if the involution on X is trivial, A(X) is just the usual Banach algebra of continuous functions on X with values in A. On the other hand, if X is a space with 2 points which are switched by the involution, A(X) is isomorphic to A'.

**7. The role of Clifford algebras.** In general, we define $S^{p,q}$ (resp. $D^{p,q}$) as the sphere (resp. the ball) of $\mathbf{R}^{p+q}$ with the involution induced by $(x_1, ..., x_p, y_1, ..., y_q) \mapsto (-x_1, ..., -x_p, y_1, ..., y_q)$ on $\mathbf{R}^{p+q}$. For p > q, the locally compact space $S^{p,0} - S^{q,0}$ is G-homeomorphic to $S^{p-q,0} \times \mathbf{R}^{q,0}$ and we have therefore the following exact sequence of Banach algebras

$$0 \longrightarrow A(S^{p-q,0})(\mathbf{R}^{q,0}) \longrightarrow A(S^{p,0}) \longrightarrow A(S^{q,0}) \longrightarrow 0$$

On the other hand, using the Clifford algebra definition of the higher K-groups, it has been proved in [3] that for any Banach algebra with unit Λ, we have natural isomorphisms

$$K^{p,0}(\Lambda) \cong K(\Lambda(\mathbf{R}^{p,0})) \cong K^p(\Lambda) \cong K(\Lambda(\mathbf{R}^{0,8k-p}))$$

for k large enough. In this formula $K^{p,q}(\Lambda)$ denotes in general the Grothendieck group of the "restriction of scalars" functor $\mathbb{P}(C^{p,q+1} \otimes \Lambda) \longrightarrow \mathbb{P}(C^{p,q} \otimes \Lambda)$. Here $C^{p,q}$ is the standard Clifford algebra of $\mathbf{R}^{p+q}$ provided with the quadratic form

$$-(x_1)^2 - ... - (x_p)^2 + (x_{p+1})^2 + ... + (x_{p+q})^2$$

and $\mathbb{P}(B)$ is the category of finitely generated projective B-modules.

In more modern and accurate homotopical terms, one may say alternatively that the homotopy fiber of the map

$$\mathcal{K}(C^{p,1} \otimes \Lambda) \longrightarrow \mathcal{K}(C^{p,0} \otimes \Lambda)$$

is also the homotopy fiber $\mathcal{F}$ of the map

$$\mathcal{K}(\Lambda) \cong \mathcal{K}(\Lambda(D^{p,0})) \longrightarrow \mathcal{K}(\Lambda(S^{p,0}))$$

One basic theorem proved in [3] § 3.4 is now the following : there is a natural homotopy equivalence $\Omega(\mathcal{K}(\Lambda(S^{p,0}))) \cong \mathcal{F} \times \Omega(\mathcal{K}(\Lambda))$ if $p \geq 3$. In other words, using Bott periodicity, we have a homotopy equivalence

$$\mathcal{K}(\Lambda(S^{p,0})) \cong \mathcal{K}(\Lambda) \times \Omega^{8k-p-1}(\mathcal{K}(\Lambda)) \cong \mathcal{K}(\Lambda) \times \Omega^{-p-1}(\mathcal{K}(\Lambda))$$

if $p \geq 3$ and $8k \geq p+1$ [with a slight abuse of notations, we write $\Omega^{-n}(\mathcal{K}(\Lambda))$ for $\Omega^{-n+8k}(\mathcal{K}(\Lambda))$, k large enough]. More precisely, we have a homotopy fibration

$$\Omega^{-p}(\mathcal{K}(\Lambda)) \longrightarrow \mathcal{K}(\Lambda(D^{p,0})) \longrightarrow \mathcal{K}(\Lambda(S^{p,0}))$$

and the first arrow is induced by the cup-product with $\eta^p$, where $\eta$ is the genererator of $\pi_1(\mathcal{K}(\mathbf{R})) \cong \mathbf{Z}/2$. It is well known that $\eta^p = 0$ when $p \geq 3$ and therefore the first arrow is null-homotopic in this case. What is proved in [3] § 3.4 is slightly more precise : there is a natural splitting $\mathcal{K}(\Lambda(S^{p,0})) \longrightarrow \mathcal{K}(\Lambda(D^{p,0}))$ from which we deduce the homotopy decomposition of $\mathcal{K}(\Lambda(S^{p,0}))$ mentionned above.

**8. Proof of the descent theorem.** We first prove by induction on p ($1 \leq p \leq 3$) that $\sigma$ induces a homotopy equivalence

$$\sigma_p : \mathcal{K}(A(S^{p,0})) \longrightarrow \mathcal{K}(A'(S^{p,0}))^{hG}$$

In this notation A'(Z) simply means the algebra of continuous functions on Z with values in A' : A'(Z) is of course the complexification of A(Z). For p = 1, this has already been shown in § 5 since $A(S^{1,0}) \cong A'$ and $(A')' \cong A' \times A'$. For p = 2, we have the following commutative diagram of homotopy fibrations (cf. § 7) :

$$\begin{array}{ccccc}
\mathcal{K}(A(S^{1,0})(\mathbf{R}^{1,0})) & \longrightarrow & \mathcal{K}(A(S^{2,0})) & \longrightarrow & \mathcal{K}(A(S^{1,0})) \\
\downarrow & & \downarrow & & \downarrow \\
\mathcal{K}(A'(S^{1,0})(\mathbf{R}^{1,0}))^{hG} & \longrightarrow & \mathcal{K}(A'(S^{2,0}))^{hG} & \longrightarrow & \mathcal{K}(A'(S^{1,0}))^{hG}
\end{array}$$

Since the two extreme vertical maps are homotopy equivalences for any Banach algebra A, it follows that the second vertical map is also a homotopy equivalence.
For p = 3, the same argument shows that $\sigma_3$ is also a homotopy equivalence.

Now, according to § 7 again, we have a commutative diagram of <u>canonically split</u> homotopy fibrations

$$\begin{array}{ccccccccc}
* & \longrightarrow & \mathcal{K}(A) & \longrightarrow & \mathcal{K}(A(S^{3,0})) & \longrightarrow & \Omega^4(\mathcal{K}(A)) & \longrightarrow & * \\
& & \downarrow & & \downarrow & & \downarrow & & \\
* & \longrightarrow & \mathcal{K}(A')^{hG} & \longrightarrow & \mathcal{K}(A'(S^{3,0}))^{hG} & \longrightarrow & \Omega^4(\mathcal{K}(A'))^{hG} & \longrightarrow & *
\end{array}$$

Since the middle map is a homotopy equivalence by our previous induction, the descent theorem follows immediately.

**9.** From the previous considerations, we can also extract a spectral sequence converging to $K^{p+q}(A(S^{3,0})) = K^{p+q}(A) \oplus K^{p+q+4}(A)$. This may be done by general methods of equivariant cohomology. More precisely, using the elementary ideas developed in [4] § 1, one can prove that the Grothendieck group $K(A(X))$ is naturally isomorphic to the group of homotopy classes of G-equivariant maps

$$X \longrightarrow \mathcal{K}(A')$$

In other words, the K-theory space $\mathcal{K}(A(X))$ may be identified - up to homotopy - with the function space of G-equivariant maps $X \longrightarrow \mathcal{K}(A')$. In particular, one has the homotopy equivalences

$$\mathcal{K}(A) \times \Omega^4(\mathcal{K}(A)) \cong \mathcal{K}(A(S^{3,0})) \cong \{S^{3,0}, \mathcal{K}(A')\}^G$$

where $\{\,,\,\}$ means function space. Since $\{S^{3,0}, \mathcal{K}(A')\}^G$ is also the space of sections of the Borel fibration

$$S^{3,0} \times_G \mathcal{K}(A') \longrightarrow S^{3,0}/G = RP^2$$

one has the following theorem :

**10. THEOREM.** *There is a spectral sequence with* $E_2$ *term* $H^p(RP_2 ; K^q(A'))$ *converging to* $K^{p+q}(A) \oplus K^{p+q+4}(A)$.

In this spectral sequence, we have of course $0 \leq p \leq 2$ and q defined mod. 4. There is at most one non zero differential which is $d_2 : H^0(RP^2 ; K^q(A')) \longrightarrow H^2(RP^2 ; K^{q-1}(A'))$ where the groups $K^*(A')$ define a local coefficient system over $RP^2$.

**11.** *Example*. If we come back to Atiyah's KR-theory of a space X, one has therefore a spectral sequence with $E_2$ term $H^p(RP^2 ; KU^q(X))$ converging to $KR^{p+q}(X) \oplus KR^{p+q+4}(X)$.

**12.** *Remark*. Let $\mathcal{KSG}(A)$ be the homotopy fiber of the map $1 - t : \mathcal{K}(A') \longrightarrow \mathcal{K}(A')$ where t is the complex conjugation. As it was essentially shown by Atiyah (see also [4] p. 178), one has a homotopy fibration

$$\mathcal{K}(A') \longrightarrow \mathcal{K}(A) \times \mathcal{K}(A'') \longrightarrow \mathcal{KSG}(A)$$

The associated homotopy exact sequence gives rise essentially to the same information as the spectral sequence above.

**13. Generalization.** The previous considerations may be used in a variety of contexts (for example to prove the Thom isomorphism in <u>real</u> K-theory, starting from the Thom isomorphism in complex K-theory) and it might be useful for the future to extract its main ideas. For this, we should consider two functors from Banach algebras to spaces, called for instance F(A) and G(A) (in our example $\mathcal{K}(A)$ and $\mathcal{K}(A')^{h\mathbf{Z}/2}$ respectively), and a natural transformation

$$\alpha : F(A) \longrightarrow G(A)$$

The claim is now the following : under suitable hypothesis, if $F(A') \cong G(A')$ by this natural transformation, then $F(A) \cong G(A)$ (isomorphisms are taken in the homotopy category). Following our previous arguments (§ 7 and 8), we see by inspection that it is enough to verify the following three conditions :

1. F and G satisfy the Mayer-Vietoris axiom (cartesian squares of Banach algebras give rise by F and G to homotopy cartesian squares).

2. $F(A(S^{1,0}))$ [resp. $G(A(S^{1,0}))$] is naturally isomorphic to $F(A')$ [resp. $G(A')$] in a way compatible with $\alpha$.

3. For p large enough, and in a way compatible with $\alpha$, F(A) [resp. G(A)] is a natural direct summand of $F(A(S^{p,0}))$ [resp. $G(A(S^{p,0}))$] through the map induced by the obvious ring homomorphism $A \longrightarrow A(S^{p,0})$.